\documentstyle[12pt]{article}
\title{Convergence of Symmetric Diffusions  on Wiener Spaces }

\author
{A. Posilicano\thanks
{Dipartimento di Scienze, Universit\`a dell'Insubria, I-22100
Como, Italy, E-mail: {\tt posilicano@uninsubria.it}}
\and  T.S. Zhang\thanks{Department of Mathematics, University of Manchester,
Manchester M139PL, UK, E-mail: {\tt tzhang@maths.man.ac.uk}}}

\begin{document}

\maketitle

\section{Introduction and Framework}

Let $(H,E,\mu )$ be an abstract Wiener space in the sense of Gross, 
i.e., $E$ is  a real separable Banach space, 
$H\subset E$ is  a continuously and densely
embedded  real separable Hilbert space, and  $\mu$ is  a centered
Gaussian  measure on $(E,{\cal B}(E))$ with covariance space $H$. The underlying  $L^p$-space is denoted by $L^p(X,\mu)$.
Define  ${\cal F}{\cal C}^\infty$ as the set of all functions $f$  of 
the form 
\begin{equation}
f(x)=\varphi (l_1(x)),\ldots ,l_n(x)),\ \varphi\in C^\infty _b(R^n), \ 
l_i\in E^{\ast}
\end{equation}
where $C^{\infty}_b(R^n)$ is the space of all bounded, infinitely differentiable
real-valued functions with bounded  partial derivatives, $E^{\ast}$ stands for the dual space of $E$.
\vskip 0.1cm

For $f\in {\cal F}{\cal C}^{\infty}$, denote by $\nabla f(x)$  the gradient of $f$
in $H$, i.e.,
\begin{equation}
<\nabla f(x),h>_H 
=\left.{{df(x+\varepsilon h)}\over {d\varepsilon}}\right|_{\varepsilon=0},
\quad \hbox{for all }\quad  h\in H,
\end{equation}
where $<\cdot,\cdot>_H$ stands for the inner product in $H$.
\vskip 0.2cm
For $p\geq 1$, denote by $D_p^1$ the completion of ${\cal F}{\cal C}^{\infty}$ under the norm
\begin{equation}
||f||_{p,1}^p=\int_E |f(x)|^p \mu(dx)+\int_E|\nabla f|_H^p(x)\mu(dx)
\end{equation}
Take $\phi \in D_{2+\varepsilon}^1$ for some $\varepsilon >0$ with  $\phi >0$, $a.s.$ and  $\int_E\phi^2(x)\mu (dx)=1$. Consider the symmetric form: 
\begin{equation}
{\cal E}_{\phi}^0(f,g)=\int_E <\nabla f(x),\nabla g(x)>_H\,\phi^2(x)\mu (dx), 
\quad f,g \in {\cal F}{\cal C}^{\infty}
\end{equation}
It is known (see e.g. [RZ1]) that $({\cal E}_{\phi}^0,D({\cal E}_{\phi}^0))$ is closable on $L^2(E,\phi^2 d\mu )$, whose closure, denoted
by $({\cal E}_{\phi},D({\cal E}_{\phi}))$, is a Dirichlet form. 
The closability of the form ${\cal E}_{\phi}$ is also equivalent to the closability of 
the gradient operator $\nabla $
in $L^2(E,\phi^2\mu)$. The action of the closure $\nabla $ on $f\in D({\cal E}_{\phi})$ will also be  denoted by $\nabla f$. 

Let $\Omega=C([0,\infty )\rightarrow E )$ be the space of all
continuous functions from $[0,\infty )$ into $E$ and let $X_t$ be the
coordinate function on $\Omega$ such that $X_t(\omega )=\omega
(t)$. The diffusion process associated with $({\cal E}_{\phi},D({\cal
E}_{\phi}))$ will be denoted by $\{\Omega, X_t, {\cal F}_t, P_x, x\in
E \}$. The Ornstein -Ulenbeck process associated with the 
Dirichlet form ${\cal E}\equiv{\cal E}_1$ will be denoted by $\{\Omega, X_t, {\cal F}_t, Q_x, x\in E \}$. Define two probability measures on $\Omega$ by

\begin{equation}
P_{\phi}(\cdot )=\int_EP_x(\cdot )\phi^2\,d\mu, \quad \quad Q_{\mu}=\int_E Q_x(\cdot )\,d\mu
\end{equation}

Let $\{\phi_n, n\geq 1\}$, $\phi_n\in D^1_{2+\varepsilon}$, be a sequence of positive (a.s. with repect to $\mu$) functions such that $\phi_n\rightarrow \phi$ in the Sobolev space $D_2^1$. Denote  the diffusion process associated with the Dirichlet form  $({\cal E}_{\phi_n},D({\cal E}_{\phi_n}))$ by $\{\Omega, X_t, {\cal F}_t, P_x^n, x\in E \}$. Define
$$
P_{\phi_n}(\cdot )=\int_EP_x(\cdot )\phi_n^2\,d\mu
$$
In this paper we will show that 
$$
\sup_{A\in {\cal F}_t}\,|P_{\phi_n}(A)-P_\phi(A)|\rightarrow 0\,,\quad t>0\,,
$$ 
i.e. $P_{\phi_n}$ converges to $P_\phi$ in total variation norm 
on ${\cal F}_t$ for any $t>0$. This
convergence is strictly stronger than weak convergence on the full
Borel $\sigma$-algebra of the path space $\Omega$.

Note that contrarily to previously known results the drifts of the 
diffusion processes can be very singular, for example where
$\phi_n=0$, $\phi=0$. The idea of our proof, following the strategy
adopted in the finite dimensional situation (see [P]), is to use
stopping times arguments to localize the diffusions in some ``good''
sets where the drifts are sufficiently regular to prove convergence by
Girsanov transform, and then to show, by capacity arguments, that the limit 
diffusion does not hit the ``bad'' sets. Such a stategy is inspired by 
Lemma 11.1.1 in [SV]. However there (see [SV], Theorem 11.1.4) 
the diffusions are then simply localized on increasing bounded balls, 
whereas we localize on the sets where the $\phi$'s are uniformly 
bounded form above and away from zero and where 
$\phi_n\rightarrow \phi$ uniformly along a subsequence.

\vskip 0.3cm
\section{Main Results}
\vskip 0.3cm
{\bf Lemma 2.1}. There exists a standard Brownian motion $B_t, t\geq 0$ taking values in  the Wiener space $E$ such that
\begin{equation}
X_t=X_0 +B_t-\int_0^t X_s ds +2\int_0^t \frac{\nabla \phi }{\phi}(X_s)ds, \quad P_{\phi}-a.e.
\end{equation}
\vskip 0.2cm
{\bf Proof.} Define
$$ B_t :=X_t-X_0+\int_0^t X_s ds -2\int_0^t \frac{\nabla \phi }{\phi}(X_s)ds$$
For $l\in E^{\ast}$, it follows by integration by parts that
$$ \int_E \frac{\partial f}{\partial l}\phi^2\,d\mu
=-\int_E \left(2\frac{\partial \phi}{\partial l}\frac{1}{\phi}-l(x) \right)f \phi^2\,d\mu$$
Thus, by Fukushima's decomposition ( see e.g. [FOT]), $l(B_t)$ is a continuous ${\cal F}_t$-martingale with
$<l(B),k(B)>_t=<l,k>_Ht$ for $l,k \in E^{\ast}$, where $<,>$ denotes
the sharp bracket of two martingales. Therefore, $B_t,t\geq 0$ is a
$E$-valued Brownain motion following  Levy's characterization. The
Lemma is proven.
\vskip 0.3cm
{\bf Lemma 2.2}. $P_{\phi}$ is absolutely continuous with respect to
$Q_{\mu}$ on ${\cal F}_t$ for $t>0$.
\vskip 0.2cm
{\bf Proof}. Without loss of generality, assume $t=1$.  Suppose first 
$$\int_0^1 \left|\frac{\nabla \phi}{\phi}\right|(X_s)ds\leq n,\quad P_{\phi}-a.e.\,.$$  
Define a new probability measure $\bar{Q}_{\mu}$ on $(\Omega, {\cal F}_1 )$ by
$$
\frac{d \bar{Q}_{\mu}}{dP_{\phi}}=
\exp\left( -2\int_0^1 \frac{\nabla \phi}{\phi}(X_s)dB_s
-2\int_0^1 \left|\frac{\nabla \phi}{\phi}\right|^2(X_s)ds\right)
$$
By the Girsanov Theorem, we see that
$$ X_t=X_0+\bar{B}_t -\int_0^t X_s ds,$$
where $\bar{B}_t=B_t+2\int_0^t \frac{\nabla \phi }{\phi}(X_s)ds$ is a Brownian motion under $ \bar{Q}_{\mu}$.
It follows from the uniqueness of the Ornstein-Ulenbeck process that
$\bar{Q}_{\mu}=Q_{\mu}$. So in this case, $Q_{\mu}$ is equivalent to 
$ P_{\phi} $. In the general case, introduce
$$
\tau_n=\inf\left\{t\geq 0, \int_0^t \left|\frac{\nabla \phi}{\phi}\right|^2(X_s)ds >n\right\}
$$
Since $$
E_{\phi}\left[\int_0^t \left|\frac{\nabla \phi}{\phi}\right|^2(X_s)ds\right]=
t\int_E |\nabla \phi|^2 (x) d\mu <\infty,\quad t\geq 0\,,$$ it follows that 
$\tau_n\rightarrow \infty$ $P_{\phi}-a.e.$ as $ n\rightarrow \infty$. Now, if $ Q_{\mu}(A)=0$, by the above discussion $P_{\phi}(A, \tau_n>1 )=0$. Therefore,
$$P_{\phi}(A)=P_{\phi}(A, \tau_n>1 )+P_{\phi}(A, \tau_n\leq 1 )=P_{\phi}(A, \tau_n\leq 1 )\leq P_{\phi}(\tau_n\leq 1 )$$
Letting $n\rightarrow \infty$ we get $P_{\phi}(A)=0$, which proves the Lemma.
\vskip 0.3cm
Let $Cap(\cdot )$, $Cap_{\phi}(\cdot )$ denote the $1$-capacities associated with the Dirichlet forms $({\cal E}, D({\cal E}))$ and $({\cal E}_{\phi}, D({\cal E}))$ (see [FOT] for details about capacities).
\vskip 0.2cm
{\bf Corollary 2.3.} $Cap(F_n)\rightarrow 0$ implies $Cap_{\phi}(F_n)\rightarrow 0$.
\vskip 0.3cm
{\bf Proof}. This follows from Lemma 2.2 and the probabilistic characterization of capacities in [FOT].

\vskip 0.3cm

In particular, we have that $ Cap_{\phi}(\phi^2\geq n )\rightarrow 0$ as $n\rightarrow\infty $ since the same is true for $Cap(\cdot )$. Moreover, we have 
\vskip 0.3cm
{\bf Lemma 2.4}. 
$$\lim_{n\rightarrow \infty}Cap_{\phi}\left(\phi\leq \frac{1}{n}\right)=0\,.$$
\vskip 0.3cm
{\bf Proof}. It was proven in [RZ1, RZ2] that the Markov uniqueness holds for the Dirichlet form $({\cal E}_{\phi},D({\cal E}_{\phi}))$, which particularly implies that 
$$\frac{1}{\phi\vee \frac{1}{n}}\in D({\cal E}_{\phi})\quad \quad \mbox{ for any } \quad n\geq 1$$
By the definition of capacity,
$$Cap_{\phi}\left(\phi\leq \frac{1}{n}\right)=
Cap_{\phi}\left(\phi\vee \frac{1}{n} \leq \frac{1}{n} \right)$$
$$=Cap_{\phi}\left(\frac{1}{\phi\vee \frac{1}{n}} \geq n \right)$$
$$\leq \frac{1}{n^2} {\cal E}_{\phi, 1}\left(\frac{1}{\phi\vee \frac{1}{n}},\frac{1}{\phi\vee \frac{1}{n}}\right)$$
$$=\frac{1}{n^2}\left[ \int_E\left(\frac{1}{\phi\vee \frac{1}{n}}\right)^2(x)\phi^2(x)\mu (dx)+\int_E \left|\nabla \frac{1}{\phi\vee \frac{1}{n}}\right|^2\phi^2(x)\mu (dx)\right]$$
$$\equiv I^n+II^n.$$
It is clear that
\begin{equation}
\lim_{n\rightarrow \infty} I^n =\lim_{n\rightarrow \infty}\frac{1}{n^2}\int_E\left(\frac{1}{\phi\vee \frac{1}{n}}\right)^2(x)\phi^2(x)\mu (dx)=0
\end{equation}
Now,
$$\nabla \frac{1}{\phi\vee \frac{1}{n}}=\left(\frac{1}{\phi\vee \frac{1}{n}}\right)^2\nabla \left(\phi\vee \frac{1}{n}\right)$$
$$=\left(\frac{1}{\phi\vee \frac{1}{n}}\right)^2\left[\nabla \phi |_{\{\phi >\frac{1}{n}\}}+\frac{1}{2}\nabla \phi |_{\{\phi =\frac{1}{n}\}}\right]$$
It follows that
$$II^n\leq \frac{1}{n^2}\int_{\left\{\phi \geq \frac{1}{n}\right\}}
\left(\frac{1}{\phi\vee \frac{1}{n}}\right)^4(x)|\nabla \phi |^2(x) \phi^2(x)\mu (dx)$$
\begin{equation}
\leq \frac{1}{n^2}\int_{\left\{\phi \geq \frac{1}{n}\right\}}\left(\frac{1}{\phi\vee \frac{1}{n}}\right)^2(x)|\nabla \phi |^2(x)\mu (dx)
\end{equation}
Observe  that 
$$ \frac{1}{n^2}|_{\{\phi \geq \frac{1}{n}\}}\left(\frac{1}{\phi\vee \frac{1}{n}}\right)^2|\nabla \phi |^2\rightarrow 0 \quad a. s.$$
as $n\rightarrow \infty$ and 
$$\frac{1}{n^2}|_{\{\phi \geq \frac{1}{n}\}}\left(\frac{1}{\phi\vee \frac{1}{n}}\right)^2|\nabla \phi |^2\leq |\nabla \phi |^2\in L^1(E,\mu )$$
By dominated convergence theorem, we have 
$$\lim_{n\rightarrow \infty}II^n=0$$
Hence, $$\lim_{n\rightarrow \infty}Cap_{\phi}\left(\phi\leq \frac{1}{n}\right)=0\,.$$
\vskip 0.3cm
Let $\phi_n, n\geq 1$ be a sequence functions in $D_2^1$ with $\phi_n>0$ $ a.s.$.
\vskip 0.3cm
{\bf Lemma 2.5}.  Assume $\phi_n\rightarrow \phi$ in the Sobolev space $D_2^1$ as $ n\rightarrow \infty$. Then there exists a decreasing sequence $\{G_m, m\geq 1\}$ of open subsets of $E$ satisfying

(i) $$\lim_{m\rightarrow \infty}Cap_{\phi}(G_m)=0\,;$$

(ii) $$\frac{1}{m}\leq \phi \leq m \quad \hbox{\rm on  $G_m^c$}\,;$$

(iii) There exists a subsequence $\{\phi_{n_k}, k\ge 1\}$ such that $$
\lim_{k \rightarrow\infty}\phi_{n_k}=\phi $$ uniformly on $G_m^c$ for
each $m$.

\vskip 0.3cm
{\bf Proof}. By Lemma 2.4, there exists a decreasing sequence
$\{G_m^{\prime}, m\geq 1\}$ of open subsets such that
$Cap_{\phi}(G_m^{\prime})\rightarrow 0$ and
$\frac{1}{m}\leq \phi \leq m$ \quad on $ (G_m^{\prime})^c$. Since
$\phi_n\rightarrow \phi $ in $D_2^1$, by Theorem 2.1.4 in [FOT] one
can find a decreasing sequence $\{G_m^{\prime\prime}, m\geq 1\}$ of
open subsets such that $Cap (G_m^{\prime\prime})\rightarrow 0$ as
$m\rightarrow \infty$ and 
$\phi_n\rightarrow \phi $ along
a subsequence, uniformly on $(G_m^{\prime\prime})^c$  for each $m$. 
By Corollary 2.3, $Cap_{\phi} (G_m^{\prime\prime})\rightarrow 0$ as $m\rightarrow \infty$. Now, set $G_m=G_m^{\prime}\cup G_m^{\prime\prime}$. Then $\{G_m, m\geq 1\}$ are the open subsets desired.

\vskip 0.3cm
Define 
\begin{equation}
\tau_m=\inf\{t>0; X_t\in G_m\}
\end{equation}
It follows from (i) in Lemma 2.5 and the probabilistic characterization of the capacity in [FOT]  that
\begin{equation}\lim_{m\rightarrow \infty}P_{\phi}(\tau_m <t )=0.
\end{equation}
\vskip 0.4cm
{\bf Theorem 2.6.} Assume that $\phi_n,\phi\in D_{2+\varepsilon}^1$ for
some positive constant $\varepsilon$ and
that $\phi_n\rightarrow \phi$ in the Sobolev space $D_{2}^1$ as $
n\rightarrow \infty$. Then $$\sup_{A\in {\cal F}_t}
\,|P_{\phi_n}(A)-P_\phi(A)|\rightarrow 0$$ for any $t>0$.
\vskip 0.3cm
{\bf Proof}. 
Define  for $m\geq 1$,
$$ \psi^{n,m}=\frac{1}{m}\vee \phi_n \wedge m \quad \quad \psi^m= \frac{1}{m}\vee \phi \wedge m$$

Denote respectively by 
$\{\Omega, X_t, {\cal F}_t, P_x^{(n,m)}, x\in E \}$ and 
$\{\Omega, X_t, {\cal F}_t, P_x^{(m)}, x\in E \}$ the diffusion
processes  associated with the Dirichlet forms  $({\cal
E}_{\psi^{n,m}},D({\cal E}_{\psi^{n,m}}))$ and $({\cal E}_{\psi^{m}},
D({\cal E}_{\psi^{m}}))$  defined as in (4) with $\phi$ replaced by 
$\psi^{n,m}$ and $\psi^{m}$. Set,
\begin{equation}
P^{(n,m)}=\int_E P_x^{(n,m)}(\cdot )(\psi^{n,m})^2d\mu, \quad P^{(m)}=\int_E P_x^{(m)}(\cdot )(\psi^{m})^2d\mu\,.
\end{equation}
For brevity let us denote by $\left\{\phi_k, k\ge 1\right\}$ the
subsequence given in Lemma 2.5(iii) and let $P^{(k,m)}$ the
corresponding measures. Note that, for $k$ suffciently large, 
\begin{equation}
P^{(k,m+2)}|_{{\cal F}_{\tau_m}}=P_{\phi_k}|_{{\cal
 F}_{\tau_m}},
 \quad \quad P^{(m+2)}|_{{\cal F}_{\tau_m}}=P_{\phi}|_{{\cal F}_{\tau_m}}\,,
\end{equation}
where ``$|_{{\cal F}_{\tau_m}}$'' stands for the restriction of the
corresponding  measure  on the $\sigma$-field ${\cal F}_{\tau_m}$,
By [ARZ], Theorem 1.3 and Remark 3.4, $P^{(n,m+2)}\sim Q_\mu$ and 
$P^{(m+2)}\sim Q_\mu$ on ${\cal F}_t$ for any $t>0$, with  
\begin{equation}
\left.\frac{dP^{(n,m+2)}}{dQ_{\mu}}\right|_{{\cal
F}_t}=L_t^{\psi^{n,m+2}}\,,\quad
\left.\frac{dQ_{\mu}}{dP^{(n,m+2)}}\right|_{{\cal
F}_t}=L_t^{1/\psi^{n,m+2}}\,,\quad
\end{equation}
and
\begin{equation}
\left.\frac{dP^{(m+2)}}{dQ_{\mu}}\right|_{{\cal
F}_t}=L_t^{\psi^{m+2}}\,,\quad
\left.\frac{dQ_{\mu}}{dP^{(m+2)}}\right|_{{\cal
F}_t}=L_t^{1/\psi^{m+2}}\,,
\end{equation}
where  
$$
L_t^\psi:=\exp\left(M_t^{\ln\psi}-\frac{1}{2}\int_0^t \left|\frac{\nabla
\psi}{\psi}\right|^2(X_s)ds 
\right)\,,
$$
and
$M_t^{\ln\psi}$ denotes the
martingale additive functional parts in the Fukushima's  decomposition
of the Dirichlet processes $\ln\psi(X_t)-\ln\psi(X_0)$, see [FOT]. 

Let us denote by $\widetilde P^{(n,m+2)}$ and $\widetilde P^{(m+2)}$
the probability measures defined by 
\begin{equation}
\left.\frac{d\widetilde P^{(n,m+2)}}{dQ_{\mu}}\right|_{{\cal
F}_t}:=L_{t\wedge\tau_m}^{\psi^{n,m+2}}\,,\quad
\left.\frac{d\widetilde P^{(m+2)}}{dQ_{\mu}}\right|_{{\cal
F}_t}:=L_{t\wedge\tau_m}^{1/\psi^{m+2}}\,.
\end{equation}
Hence 
\begin{equation}
\widetilde P^{(k,m+2)}|_{{\cal F}_{\tau_m}}=P^{(k,m+2)}|_{{\cal F}_{\tau_m}}
=P_{\phi_k}|_{{\cal F}_{\tau_m}},
\end{equation}
\begin{equation}
\widetilde P^{(m+2)}|_{{\cal F}_{\tau_m}}=P^{(m+2)}|_{{\cal F}_{\tau_m}}=
P_{\phi}|_{{\cal F}_{\tau_m}}\,,
\end{equation}
\begin{equation}
\widetilde P^{(n,m+2)}\sim \widetilde P^{(m+2)}\quad
\hbox{\rm on ${\cal F}_t$ for any $t>0$,}
\end{equation}
and 
\begin{eqnarray}
&\ &\left.\frac{dP^{(n,m+2)}}{dP^{(m+2)}}\right|_{{\cal
F}_t}=\exp\left(M_{t\wedge\tau_m}^{\ln\psi^{,m+2}}-M_{t\wedge\tau_m}^{-\ln\psi^{m+2}}\right.\nonumber\\
&\ &\left.-\frac{1}{2}\int_0^{t\wedge\tau_m} \left(\left|\frac{\nabla
\psi^{n,m+2}}{\psi^{n,m+2}}\right|^2-\left|\frac{\nabla
\psi^{m+2}}{\psi^{m+2}}\right|^2\right)(X_s)ds \right)\,.
\end{eqnarray}
Proceeding as in the proof of Lemma 3.1 in [DP] (a sort of versions
in total variation norm of Lemma 11.1.1 in [SV]) one easily gets
\begin{equation}
\sup_{A\in {\cal F}_t}|P_{\phi_k}(A)-P_\phi(A)|\le 
3\sup_{A\in {\cal F}_t}|\widetilde P^{(k,m+2)}(A)-\widetilde
P^{(m+2)}(A)|+4P_\phi(\tau_m<t)
\end{equation}
and, by the Csizl\'ar-Kullback inequality and (19), one obtains
\begin{eqnarray*}
&\ &\left(\sup_{A\in {\cal F}_t}|\widetilde P^{(k,m+2)}(A)
-\widetilde P^{(m+2)}(A)|\,\right)^2\le
2E_{\widetilde
P^{(k,m+2)}}\left[\ln\left.\frac{d\widetilde P^{(k,m+2)}}{d\widetilde 
P^{(m+2)}}\right|_{{\cal
F}_t}\right]\\
&=&E_{\widetilde P^{(k,m+2)}}\left[\int_0^{t\wedge\tau_m} \left(\left|\frac{\nabla
\psi^{k,m+2}}{\psi^{k,m+2}}\right|^2-\left|\frac{\nabla
\psi^{m+2}}{\psi^{m+2}}\right|^2\right)(X_s)ds\right]\\
&=&\int_0^tE_{P_{\phi_k}}\left[1_{\{\tau_m>s\}} \left(\left|\frac{\nabla
\psi^{k,m+2}}{\psi^{k,m+2}}\right|^2-\left|\frac{\nabla
\psi^{m+2}}{\psi^{m+2}}\right|^2\right)(X_s)\right]ds\\
&\le& t\|\nabla
\phi_k-\nabla\phi\|^2_{L^2(E,\mu)}+tm^2\|\nabla\phi\|^2_{L^2(E,\mu)}
\|\phi_k-\phi\|^2_{L^\infty(G_m^c,\mu)}\,.
\end{eqnarray*}
Thus by Lemma 2.5(iii), (10) and (20) 
$\sup_{A\in {\cal F}_t}|P_{\phi_n}(A)-P_\phi(A)|$
converges along a subsequence. Suppose now that the whole sequence
does not converge. Then there exists a subsequence $\{P_{\phi_{n_j}},
j\ge 1\}$ such that $$\sup_{A\in {\cal
F}_t}|P_{\phi_{n_j}}(A)-P_\phi(A)|>\varepsilon\,,\qquad
\hbox{\rm for any $j$.}$$ Since $\phi_{n_j}\to\phi$ in $D^1_2$, we get a
contradiction. This completes the proof.

\end{document}